\documentclass[12pt]{article}
  \usepackage{inputenc}
  \usepackage[english]{babel}
  \usepackage{amsfonts,amssymb,amsmath,amstext,amsbsy,amsopn,amscd,graphicx}
  \usepackage{graphics}
  \usepackage{amsthm}

 \newtheorem{thm}{Theorem}[section]
  \newtheorem{st}{Statement}[section]

  \newtheorem{dfn}{Definition}[section]
  \newtheorem{rk}{Remark}[section]

 \def\R{{\Bbb R}}
 \def\Z{{\Bbb Z}}

\title{On Virtual Crossing Number Estimates For Virtual Links}

\author{Denis Afanasiev\footnote{Corresponding author}\,  and Vassily Manturov}

\date{}

\begin{document}

\newcommand{\skcrro}{\raisebox{-0.25\height}{\includegraphics[width=0.5cm]{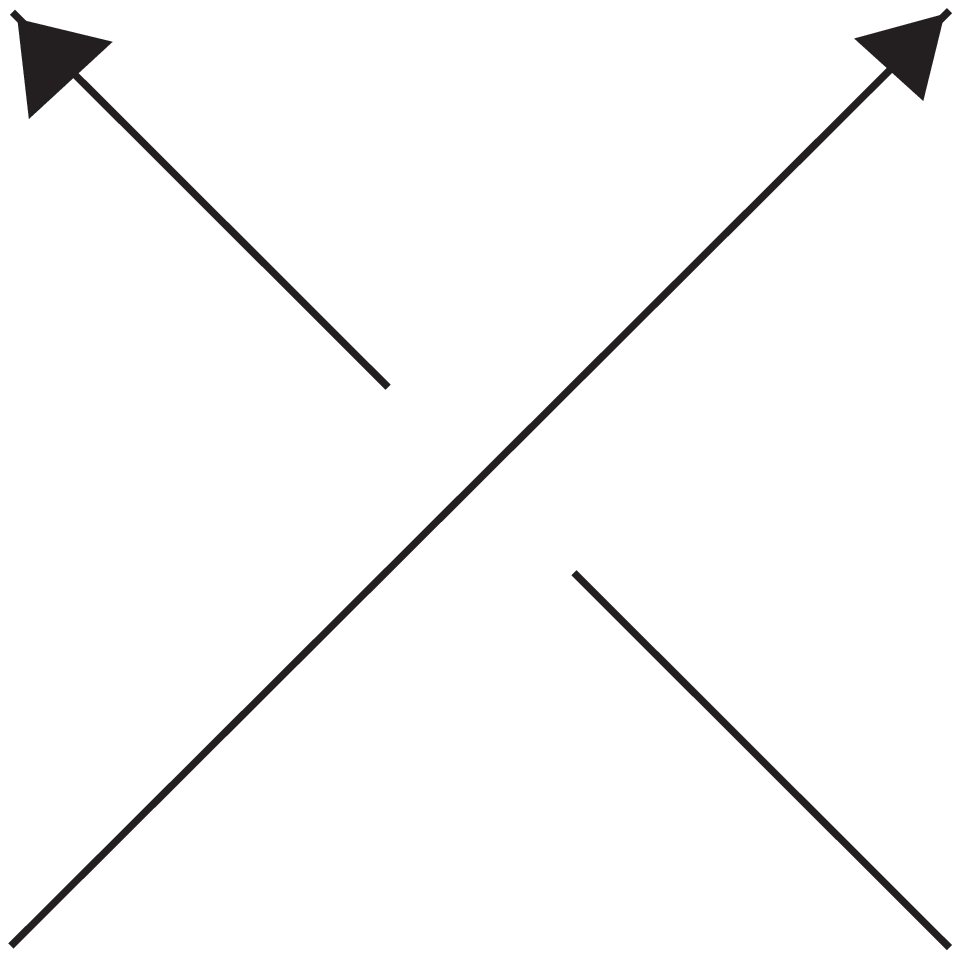}}}
\newcommand{\skcrlo}{\raisebox{-0.25\height}{\includegraphics[width=0.5cm]{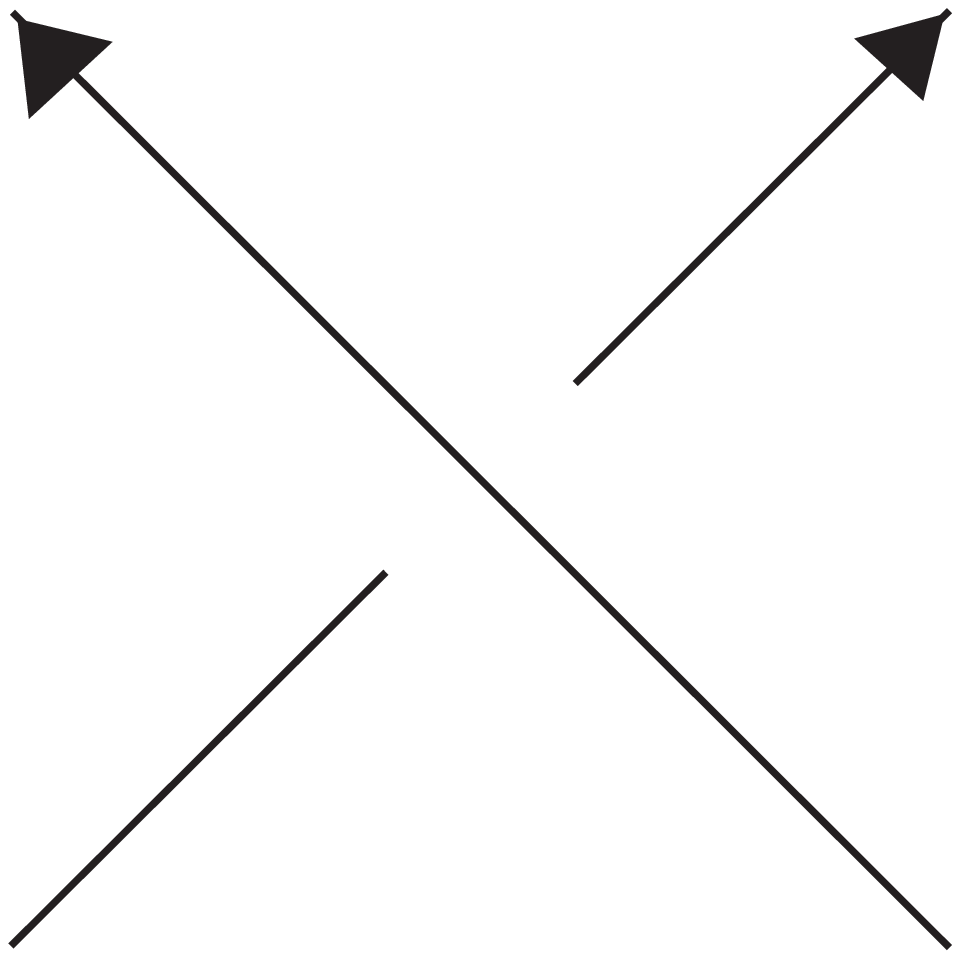}}}

\maketitle

\begin{abstract}

We address the question of detecting minimal virtual diagrams with
respect to the number of virtual crossings. This problem is closely
connected to the problem of detecting the minimal number of
additional intersection points for a generic immersion of a singular
link in $R^{2}$.

We tackle this problem by the so-called $\xi$-polynomial whose
leading (lowest) degree naturally estimates the virtual crossing
number.

Several sufficient conditions for minimality together with infinite
series of new examples are given.

We also state several open questions about $M$-diagrams, which are
minimal according to our sufficient conditions.

\end{abstract}

\section{Introduction}

One of the most important problems in the classification of knots
and links is the problem of detecting the minimal crossing number.
One of the striking achievements in that direction is the celebrated
Kauffman-Murasugi theorem \cite{Mur,Thi}, for further
generalizations see \cite{MyBook} establishing the minimality of
reduced alternating diagrams for virtual links. This theorem
(conjectured by Tait in late 19-th century) was solved in late
1980-s by using a newborn polynomial: the Kauffman bracket version
\cite{KauffmanBracket} of the Jones polynomial \cite{Jones}.

In mid-90-s, Louis Kauffman introduced \cite{KaV}  a natural
generalization of knots and links called virtual knot theory. This
generalization can be treated both topologically (as knots in
thickened 2-surfaces up to homotopy and
stabilization/destabilization) and combinatorially (via a
generalization of planar diagrams with a new crossing type ---
called virtual --- allowed). Virtual knot theory played a crucial
role in solving some problems about classical knots (see, e.g.,
\cite{GPV}) and for understanding of some notions in classical knot
theory, and stating new ones (see, e.g.\cite{FKM}).

One of the most natural problems in virtual knot theory is to
estimate classical and virtual crossing numbers. For estimating
classical crossing number, a lot of results were obtained (see,
e.g., \cite{Minim,Miyazawa,MyBook} and references therein) by using
some generalization of the Kauffman-Murasugi theorem etc.

The problem this present paper is devoted to is to estimate from
below the virtual crossing number: a partial case of this theorem
detects non-classicality of a link (in this direction we note the
paper \cite{DK}). Note that this problem is closely connected to the
problem of estimating minimal number of crossings for projecting a
given singular link to a plane.

We attack the virtual crossing number by using the $\xi$-polynomial
introduced independently by several authors (see
\cite{KR},\cite{Saw},\cite{SW},\cite{MaXi}, for the proof of their
coincidence, see \cite{BF}, for detecting non-classicality see also
\cite{Tep}), and we use the definition by the second named author of
the present paper: this definition includes a ``counting variable''
for virtual crossings, and it is natural that the degrees of the
polynomial invariant in this counting variable estimates the virtual
crossing number from below.

The paper is organized as follows.

The next section contains all basic notions and the main Theorem
(Theorem 1) on virtual crossing estimates.

In section 3, we prove some theorems delivering sufficient
conditions for minimality and present several series of examples of
minimal diagrams.

We conclude our paper by section 4 with a list of unsolved problems.

\subsection{Acknowledgements}

The authors are grateful to L.H.Kauffman, A.T.Fomenko and
V.A.Vassiliev for their attention to this work and fruitful
consultations.

\section{Basic Notions and Constructions}

Virtual knot theory was invented by Kauffman around 1996,
\cite{KaV}.

 \begin {dfn}
 {\upshape
A {\em virtual diagram} is a $4$-valent diagram in $\R^2$ where each
crossing is either endowed with a classical crossing structure (with
a choice for underpass and overpass specified) or just said to be
virtual and marked by a circle.}
 \end {dfn}

\begin{figure}
 \centering\includegraphics[width=200pt]{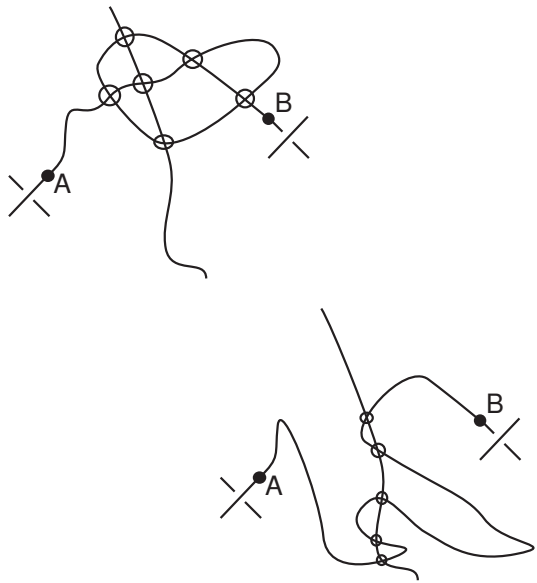}
  \caption{The detour move}
  \label{detour}
 \end{figure}

 \begin {dfn}{\upshape
A {\em virtual link} is an equivalence class of virtual link diagram
modulo generalised Reidemeister moves. The latter consist of usual
Reidemeister moves referring to classical crossings and the {\em
detour move} that replaces one arc containing only virtual
intersections and self-intersection by another arc of such sort in
any other place of the plane, see Fig.~\ref{detour}.}
 \end {dfn}

Throughout the paper, we deal only with oriented links. By an {\em
arc} of a planar diagram we mean a connected component of the set,
obtained from the diagram by deleting all virtual crossings (at
classical crossing the undercrossing pair of edges of the diagram is
thought to be disjoint as it is usually illustrated). We say that
two arcs $a, a'$ belong to the same {\em long arc} if there exists a
sequence of arcs $a = a_{0},\dots,  a_{n} = a'$ and virtual
crossings $c_{1},\dots,c_{n}$ such that for $i = 1,\dots, n$ the
arcs $a_{i},a_{i+1}$ are incident to $c_i$ from opposite sides.  A
planar virtual diagram is said to be {\em proper}, if it has no
cyclic long arcs. It is easy to show that equivalent proper diagrams
can be transformed to each other by using generalized Reidemeister
moves in the class of proper diagrams. This can be done by
adding/removing a curl (by using the first classical Reidemeister
move) when necessary. For each proper diagram, the number of long
arcs equals the number of classical crossings. In the sequel, we
deal only with proper diagrams.

Now, let us construct the $\xi$-polynomial of virtual links (and its
unnormalized version called $\zeta$).

Let $D$ be a proper diagram of a virtual link $D$ with $n$ classical
crossings. Let us construct an $n\times n$-matrix $A(D)$ with
elements from $\Z[t,s,t^{-1},s^{-1}]$ as follows.

First,we enumerate all classical crossings of $D$ by integer numbers
from $1$ to $n$ and associate with each crossing the ematnaing long
arc. Each long arc starts with a (short) arc. Let us associate the
label 1 with latter arc. All other arcs of the long arc will be
marked by exponents $s^{k}, k \in \Z$, as follows. While passing
through the virtual crossing, we multiply the label by $s$ if we
pass from the left to the right or by $s^{-1}$ otherwise.

Since the diagram is proper, our labeling is well defined. Consider
a classical crossing  $v_i$ with number $i$. It is incident to some
three arcs $p, q, r$, belonging to long arcs with numbers $i, j, k$,
whence number $j$ belongs to the arc passing through $v_i$. Denote
the power of $s$ of the label corresponding to $q$ by $a_{ij}$, and
that of the label corresponding to $r$ by $a_{ik}$.

Let us define the $i$-th row of the matrix $A(D)$ as the sum of the
following three rows $y_{1}, y_{2}, y_{3}$ of length $n$. Each of
these rows has only one nonzero element. The $i$-th element of the
row $y_{1}$ equals 1. If the crossing is positive, we set

\begin{equation}
y_{2k}=-s^{a_{ik}}t,y_{3j}=(t-1)s^{a_{ij}};
\end{equation}

otherwise we put

\begin{equation}
y_{2k}=-s^{a_{ik}}t^{-1},y_{3j}=(t^{-1}-1)s^{a_{ij}}
\end{equation}

We set $\zeta(D):=det\,A(D)$. In \cite{MaXi} it is proved that if
two virtual diagrams $D,D'$ are equivalent then $\zeta(D)=
t^{l}\zeta(D')$ for some integer $l$.

Denote by $deg f$ and $mdeg f$ the leading (lowest) degree of
monomials $f\in \mathbb{Z}[t,t^{-1}][s,s^{-1}]$ with respect to
variable $s$; if $f=0$, we set $deg f=-\infty$, $mdeg f=+\infty$.

\begin{thm}\label{theorem1} Let $k$ be the number of virtual crossings
of a virtual diagram $D$. Then $deg\,\zeta(D)\leqslant k$,
$mdeg\,\zeta(D)\geqslant -k$.
\end{thm}

\begin{rk} $deg\,\zeta(D)$ and $mdeg\,\zeta(D)$ are virtual link invariants.
\end{rk}

\begin{proof}[Proof] Let $v_1,...,v_n$ be an arbitrarily enumerated
classical crossings of a virtual diagram $D$,
$\gamma^1,...,\gamma^n$ be long arcs of $D$ enumerated according to
the rule: $\gamma^i$ is emanating from $v_i$. Here $k_i$ denotes
number of virtual crossing that belong to $\gamma^i$,
$\gamma^i_0,...,\gamma^i_{k_i}$ denote arcs on $\gamma^i$ enumerated
in their sequence on the long arc. Now, using this notation, we can
get another definition of Alexander-like matrix $A=A(D)$.

{\upshape The polynomial $[v_i:\gamma_{\mu}^j]\in
\mathbb{Z}[t,t^{-1}]$,
 defined in following way, is called an {\em incidence coefficient} of classical
crossing $v_i$ and arc $\gamma_{\mu}^j$. If $\gamma_{\mu}^j$ and
$v_i$ are not incident, we set $[v_i:\gamma_{\mu}^j]= 0$.
 If they are incident, we have the following
 options:   the arc (1) is emanating from the crossing, (2) or passing
through it, (3) or coming into it, or simultaneosly realize some of
conditions (1)-(3); accordingly,
 $[v_i:\gamma_{\mu}^j]$ is defined as a sum of some of three polynomials:

 (1) $1$;

 (2) $t^{sgn\, v_i}-1$;

 (3) $-t^{sgn\, v_i}$;

\noindent where $sgn\,v_i$ denotes a {\em local writhe number} of
crossing $v_i$ (Recall that the {\em local writhe number} of a
classical crossing is $+1$ for $\skcrro$ and $-1$ for $\skcrlo$. The
{\em writhe number} of a (virtual) diagram is the sum of local
writhe numbers over all classical crossings). }

 Thus, we set
$A_{ij}:=\sum_{\mu=0}^{k_j}\,[v_i:\gamma_{\mu}^j]s^{deg\,\gamma_{\mu}^j
}$, where $deg\,\gamma_{\mu}^j$ denotes the power of variable $s$ on
 arc $\gamma_{\mu}^j$.

If we have a long arc passing through a virtual crossing, then this
passing may be either {\em increasing} or {\em decreasing} (with
respect to the power of $s$ on the arcs of the long arc separated by
the virtual crossing in question. More precisely, we say that a
virtual crossing $v$ is {\em increasing} on a long arc $\gamma$
between two consecutive arcs $\gamma_{i}$ and $\gamma_{i+1}$ if the
other long arc passing through $v$ goes from the left to the right
with respect to the orientation $\gamma_{i}$ and $\gamma_{i+1}$. If
the other long arc goes from the left to the right, we call this arc
{\em decreasing}.

We note that each virtual crossing is increasing for exactly one
long arc and decreasing for exactly one arc. These two arcs may
coincide if a long arc passes twice through a virtual crossing.

Let $p_j$ and $q_j$, $j=1,...,n$, be numbers of increasing and
decreasing virtual crossings on $\gamma^j$, respectively. Then
$max_{0\leqslant\mu\leqslant k_j}deg\,\gamma_{\mu}^j\leqslant p_j$
and $min_{0\leqslant\mu \leqslant k_j}deg\,\gamma_{\mu}^j\geqslant
-q_j$. Consequently, $deg\,A_{ij}=
deg(\sum_{\mu=0}^{k_j}\,[v_i:\gamma_{\mu}^j]s^{deg\,\gamma_{\mu}^j
})\leqslant max_{0\leqslant\mu\leqslant k_j}deg\,\gamma_{\mu}^j
$$\leqslant p_j$ and $mdeg\,A_{ij}=$ $mdeg(\sum_{\mu=0}^{k_j}$
$[v_i:\gamma_{\mu}^j]s^{deg\,\gamma_{\mu}^j })$ $\geqslant
min_{0\leqslant\mu\leqslant k_j}deg\,\gamma_{\mu}^j$ $\geqslant
-q_j$. So, we have inequalities:

\begin{equation}\label{eneq3}
deg\,A_{ij}\leqslant max_{0\leqslant\mu\leqslant
k_j}deg\,\gamma_{\mu}^j  \leqslant p_j, mdeg\,A_{ij}\geqslant -q_j.
\end{equation}

Hence, $deg\,\zeta(D)= deg({\sum}_{\sigma\in
S_n}(-1)^{\sigma}A_{{\sigma}_1 1}A_{{\sigma}_2 2}...A_{{\sigma}_n
n})$ $\leqslant max_{\sigma\in S_n}(deg\,A_{{\sigma}_1
1}+...+deg\,A_{{\sigma}_n n})$ $\leqslant p_1+...+p_n$. Analogously,
$mdeg\,\zeta(D)\geqslant$  $ min_{\sigma\in S_n}(mdeg\,A_{{\sigma}_1
1}+...+mdeg\,A_{{\sigma}_n n})\geqslant$ $-q_1-...-q_n$.

Since every virtual crossing is increasing for exactly one long arc
and decreasing for exactly one long arc (possibly, the same), we get
 \begin{equation}\label{eq4}
   p_1+...+p_n=k, q_1+...+q_n=k.
 \end{equation}

 Thus, $deg\,\zeta(D)\leqslant p_1+...+p_n=k$ and $mdeg\,\zeta(D)\geqslant -q_1-...-q_n =
 -k$.
 The theorem is
proved.\end{proof}

\section{Sufficient conditions for minimality}
Assume for a proper virtual diagram $D$ we have
$deg\,\zeta(D)=k(D)$. Then for every virtual diagram $D'$ which
equivalent to $D$ we have $\zeta(D')= t^{l}\zeta(D)$ $\Rightarrow$
$k(D)=deg\,\zeta(D)=deg\,\zeta(D')\leqslant k(D')$ by Theorem
\ref{theorem1}. Thus, $D$ has minimal possible number of virtual
crossings.

So, sufficient conditions for {equality} $deg\,\zeta(D)=k(D)$ to
hold are also sufficient for minimality of $D$ {with respect to} the
number of virtual crossings.

Now, our aim is to find necessary and sufficient conditions, when
the equation $deg\,\zeta(D)=k(D)$ holds. Let $A^l_{ij}$,
$l\in\mathbb{Z}$, be the coefficients of the polynomial $A_{ij}\in
\mathbb{Z}[t,t^{-1}][s,s^{-1}]$ defined as follows.
 $A_{ij}= \sum_{l=-\infty}^{+\infty}A^l_{ij}s^l.$
We shall use the notation $(f)_k$ for the $k$-th coefficient of
$f\in R[x,x^{-1}]$. Then, $(\zeta(D))_k=(det\,A(D))_k$
$=({\sum}_{\sigma\in S_n}(-1)^{\sigma} A_{{\sigma}_1 1}$ $
A_{{\sigma}_2 2}...A_{{\sigma}_n n})_k=$ ${\sum}_{\sigma\in
S_n}(-1)^{\sigma}$ $\sum_{l_1,...,l_n:\,l_1+...+l_n=k}$
$A^{l_1}_{{\sigma}_1 1}...A^{l_n}_{{\sigma}_n n}=$
$\sum_{l_1,...,l_n:\,l_1+...+l_n=k}$ ${\sum}_{\sigma\in
S_n}(-1)^{\sigma}$ $A^{l_1}_{{\sigma}_1 1}...A^{l_n}_{{\sigma}_n
n}=$ ${\sum}_{\sigma\in S_n}(-1)^{\sigma}A^{p_1}_{{\sigma}_1
1}...A^{p_n}_{{\sigma}_n n}$, because $A^l_{ij}=0$ for $l>p_j$ by
(\ref{eneq3}) and $p_1+...+p_n=k$ by (\ref{eq4}). So,
$(\zeta(D))_k=det\,T$, where $T_{ij}:=A^{p_j}_{ij}$.

Thus, to reach the maximal possible degree (minimal possible degree)
of the $\xi$-polynomial, we need $T$ to be non-degenerate, which, in
turn, is possible only when each long arc contributes its maximal
possible (minimal possible) degree of $s$ to a prefixed summand of
the determinant of $T$. More exactly, we have the following.

\begin{equation}\label{eq5}
deg\,\zeta(D)=k(D) \Leftrightarrow det\, T(D)\neq 0.
\end{equation}

 The matrix $T=T(D)$ {satisfies} following

\begin{st}\label{st1}
 1) $det\, T(D)\neq 0$ $\Longrightarrow$ $\forall j \
max_{0\leqslant\mu\leqslant k_j}deg\,\gamma_{\mu}^j = p_j$;

2) $\forall j \ max_{0\leqslant\mu\leqslant k_j}deg\,\gamma_{\mu}^j
= p_j$ $\Longrightarrow$ $T_{ij}(D)=[v_i:{\gamma}^j_{p_j}]$.
\end{st}

\begin{proof}[Proof]
(1) Let us suppose that $\exists\, j$ $max_{0\leqslant\mu\leqslant
k_j}deg\,\gamma_{\mu}^j < p_j$. Then $deg\,A_{ij}=
deg(\sum_{\mu=0}^{k_j}\,[v_i:\gamma_{\mu}^j]s^{deg\,\gamma_{\mu}^j
})\leqslant max_{0\leqslant\mu\leqslant k_j}deg\,\gamma_{\mu}^j $ $<
p_j$ $\Longrightarrow$ $\forall i\ deg\,A_{ij} < p_j$
$\Longrightarrow$ $\forall i\ T_{ij}=A^{p_j}_{ij} = 0$
$\Longrightarrow$  $det\, T(D)= 0$.

(2) $max_{0\leqslant\mu\leqslant k_j}deg\,\gamma_{\mu}^j = p_j$
$\Longrightarrow$ $\forall j\ deg \,\gamma_{p_j}^j=p_j$ and $deg
\,\gamma_{\mu}^j< p_j$ for $\mu\neq p_j$ $\Longrightarrow$ $T_{ij}=
A^{p_j}_{ij} =$
$(\sum_{\mu=0}^{k_j}\,[v_i:\gamma_{\mu}^j]s^{deg\,\gamma_{\mu}^j})_{p_j}=[v_i:{\gamma}^j_{p_j}]$.
\end{proof}

Statement \ref{st1} {naturally leads to} the following

\begin{dfn}
{\upshape
 (1) An arc $\gamma_{p_j}^j$, $j=1,...,n$, $deg\,\gamma_{p_j}^j=p_j$,
 we called {\em critical arc}. It is an arc which has a greatest
 possible degree.

 (2) Virtual diagram $D$ is called {\em special}, if every long arc
 of $D$ contains a critical arc, i.e. $\forall j \ max_{0\leqslant\mu\leqslant k_j}deg\,\gamma_{\mu}^j
= p_j$. In other words, increasing virtual crossings on each long
arc are located in the beginning of this long arc.

 (3) For a special diagram $D$, the $n\times n$-matrix $T=T(D)
 \in \mathbb{Z}[t,t^{-1}]$,
 $T_{ij}=[v_i:{\gamma}^j_{p_j}]$, is called the {\em $T$-matrix of special diagram}
 $D$. It is composed of incidence
 coefficients of classical crossings and critical arcs.

 (4) We call a special diagram $D$ with a nondegenerate $T$-matrix  a {\em
 T-diagram}.
 }
\end{dfn}

From the definition above, statement \ref{st1} and relation
(\ref{eq5}), we easily conclude

\begin{st}\label{st2}
 A proper virtual diagram is $T$-diagram if and only if
 $deg\,\zeta(D)=k(D)$. Besides, for a special virtual diagram $D$ the
 $k$-th coefficient of polynomial $\zeta=\zeta(D)$ is equal to
 determinant of T-matrix of this special diagram.
\end{st}

\begin{figure}
\centering\includegraphics[scale = 0.8]{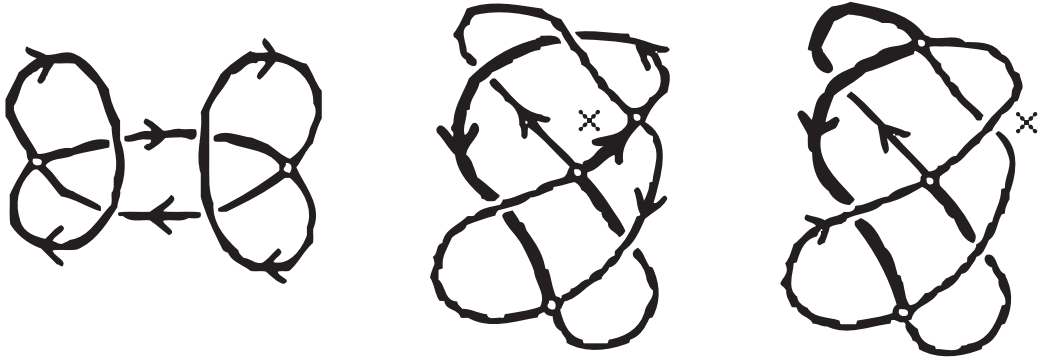}
  \caption{Special diagrams with $deg\,\zeta<k$}
  \label{first}

\end{figure}

Hence, if $D$ is not special or it is special, but there is a
critical arc which is not incident with some of classical crossing,
or there is a classical crossing which is not incident some of
critical arc, then $deg\,\zeta(D)<k$. These options are illustrated
in Figure \ref{first}.

Besides, by using Statement \ref{st2}, we formulate an easier
sufficient condition of minimality than come from Statement
\ref{st2}. Let us introduce following

\begin{dfn}
{\upshape
 For a special diagram $D$, we call the $n\times n$-matrix $M=M(D)$,
which indicates when the $i$-th classical crossing and the $j$-th
critical arc (i.e. $\gamma_{p_j}^j$) of  $D$ are incident, the {\em
incidence matrix of $D$}. }
\end{dfn}

So, if $T_{ij}=[v_i:{\gamma}^j_{p_j}]\neq 0$, then $M_{ij}= 1$. But
the inverse statement is not true at all. {For example} if
${\gamma}^j_{p_j}$ is emanating from $v_i$, passing through it and
coming into it, then $[v_i:{\gamma}^j_{p_j}]=0$.

\begin{dfn}
{\upshape
 We say that a classical crossing $v_i$ of $D$ is {\em
cyclic}, if some arc of $D$ is emanating from $v_i$, passing through
it and coming into it. }
\end{dfn}

Thus, if a special diagram $D$ does not have cyclic crossings,  then
$T_{ij}=[v_i:{\gamma}^j_{p_j}]\neq 0$ $\Leftrightarrow$ $M_{ij}\neq
0$. Consequently, if the permanent\footnote{$per\,
M:={\sum}_{\sigma\in S_n} M_{{\sigma}_1 1}M_{{\sigma}_2
2}...M_{{\sigma}_n n}$}\, of matrix $M$ is equal to $1$, in other
words, there is unique pairing of classical crossings and critical
arcs, then $det\,T(D)\neq 0$, consequently,  $D$ is minimal by
Statement \ref{st2}. Thus, we get

\begin{thm}\label{theorem2}
If the permanent of the incidence matrix of a special diagram, which
does not have cyclic classical crossings, is equal to $1$ (i.e.
there is unique pairing of classical crossings and critical arcs),
then this virtual diagram is minimal with respect to the number of
virtual crossings.
\end{thm}

\begin{figure}
  \centering\includegraphics[scale = 0.8]{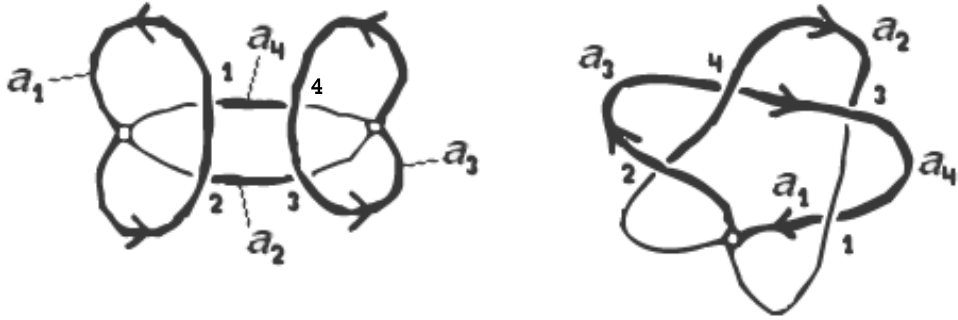}
  \caption{Special diagrams with $per M > 1$}
\label{second}
\end{figure}

In Figure \ref{second} we illustrate that the sufficient condition
for minimality, obtained from Statement \ref{st2}, is much more
difficult than the sufficient condition from Theorem \ref{theorem2}.

\noindent For the left diagram in Figure \ref{second}, we have

 $det\,T=
\begin{vmatrix} t-1 & 0 & 0 & -t\\
                {t^{-1}-1} & 1 & 0 & 0\\
                0 & {-t} & {t-1} & 0\\
                0 & 0 & {t^{-1}-1} & 1\\
 \end{vmatrix}
 =(t-1)^2-(-t)^2(t^{-1}-1)^2=0$.

\noindent Consequently, $deg\,\zeta < 2$ by Statement \ref{st2}. So,
we can not conclude that the diagram is minimal. For the right
diagram in Figure 3 we have

$det\,T=
\begin{vmatrix} 1 & 0 & 0 & -t^{-1}\\
                0 & 1 & {t^{-1}-1} & 0\\
                0 & {-t^{-1}} & 0 & {t^{-1}-1}\\
                0 & {t^{-1}-1}& {-t^{-1}} & 1\\
 \end{vmatrix}
 =
 (t^{-1}-1)t^{-1}-(t^{-1}-1)(-t^{-1}-(t^{-1}-1)^2)=(t^{-1}-1)(t^{-2}+1)=t^{-3}-t^{-2}+t^{-1}-1\neq
 0$. Thus, the diagram is minimal.

 Theorem \ref{theorem2} naturally leads to
 \begin{dfn}
 {\upshape
 A special diagram, which does not have cyclic classical crossings
 and has unique pairing of classical crossings
and critical arcs, we call an {\em $M$-diagram}.
 }
 \end{dfn}

 For detecting nonequivalence of $M$-diagrams it is convenient to
 use the
 following two invariants. Let virtual link $L$ have an $M$-diagram $D$.
Then by Statement \ref{st2} and the property of $\zeta$-polynomial
we have $det\, T(D')=t^l det\, T(D)$, where $D$ and $D'$ are
$M$-diagrams of $L$. But $det\,T(D)=det\,M(D)\, t_{1 j_1}...t_{n
j_n}$ (for some $j_1,...,j_n$) can be rewritten in the form
$\varepsilon\, t^{\alpha}(t-1)^{\beta}$, where $\varepsilon=\pm 1$,
$\beta\geqslant 0$, $\alpha,\beta\in\mathbb{Z}$, because
$t_{ij}=1,\, t^{\pm 1}-1$ or $-t^{\pm 1}$  and
$\,t^{-1}-1=(-t^{-1})(t-1)$. Note that numbers $\varepsilon,\alpha$
and $\beta$ quite defined for a product of polynomials $\pm 1,\,
t^{\pm 1}-1,\, -t^{\pm
  1}$. Thus, numbers $\varepsilon$ and $\beta$ are invariants for
  virtual link $L$.

We can formulate the  definition of the invariants $\varepsilon$ and
$\beta$ in a more geometric way. $M$-diagram $D$ has unique pairing
of classical crossings and critical arcs. In the sense of this
pairing, $\beta(D)$ is the number of classical crossings paired with
critical arcs which are passing through them.
$\varepsilon(D)=det\,M(D)\,(-1)^{x(D)+y(D)}$, where $x(D)$ is the
number of negative classical crossings (i.e. their local writhe
number is equal to $-1$) which paired with critical arcs that are
passing through they, $y(D)$ is the number of classical crossings
paired with critical arcs which are coming into they, $det\,M(D)=\pm
1$ characterizes oddness of the pairing of the $M$-diagram.

\begin{figure}
\centering\includegraphics[scale = 0.8]{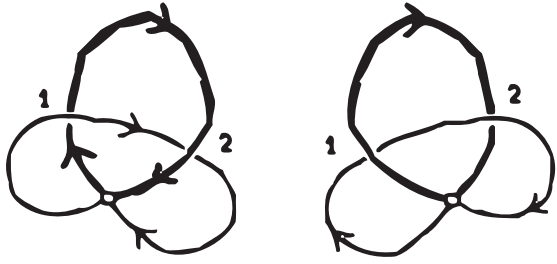}
  \caption{$M$-knots $T$ and $\widetilde{T}$}
\label{fig3}
\end{figure}

\begin{dfn}
{\upshape We say that virtual link is {\em $M$-link}, if can be
represented by an $M$-diagram.
 }
\end{dfn}
By using the $\varepsilon$-invariant, we prove that all $M$-knots
shown in Figure \ref{fig3} are different.

\noindent Indeed, we have $M(T)=\begin{pmatrix} 1 & 1\\
                1  & 0\\
                \end{pmatrix}$, $det\, M=-1$,
                  $x=0$, $y=1$ $\Rightarrow$
                $\varepsilon(T)=det\,M\,(-1)^x\,(-1)^y=1;$

\begin{figure}
\centering\includegraphics[scale = 0.8]{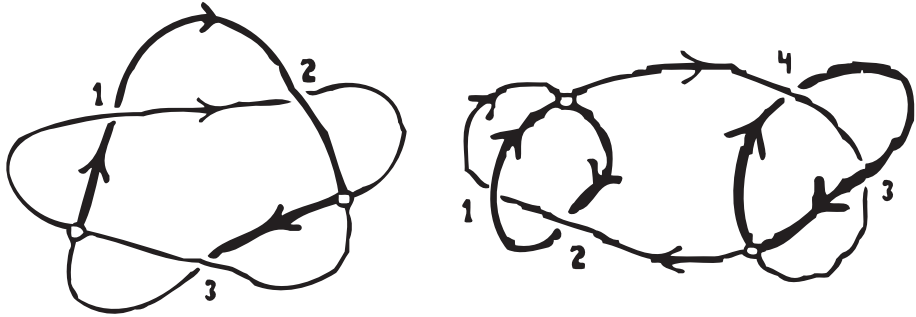}
  \caption{$M$-knots $Q$ and $T^{\prime}_2$}
\label{fig4}
\end{figure}

                $\varepsilon(\widetilde{T})=det\,\begin{pmatrix} 1 & 0\\
                1  & 1\\
                \end{pmatrix}\, (-1)^{1}(-1)^{0}=-1$.
The virtual knots shown in Figure \ref{fig4} are distinguished by
the $\beta$-invariant.

\noindent Indeed, $\beta(Q)=1$, $\beta(T^{\prime}_2)=2$, but
$\varepsilon$-invariant does not distinguish between these knots, because we have
 $\varepsilon(Q)= det\,\begin{pmatrix} 1 & 0& 1\\
                1  & 0&0\\
                0&1&0\\
                \end{pmatrix}$$\,(-1)^{0}\,(-1)^{2}=1$ and
                $\varepsilon(T^{\prime}_2)={\varepsilon(T)}^2=1$.
                Note that the equation
                $\varepsilon(T^{\prime}_2)={\varepsilon(T)}^2$ holds
                by the following construction of {\em special sum} of
                special diagrams.

For construction new $T$- and $M$-diagrams it is convenient to use
following operation of {\em connected sum of special diagrams}.
There are two types of arcs, which do not critical: {\em
pre-critical}, which precede a critical arc on a long arc, and {\em
post-critical}, which follow a critical arc.

A connected sum of special diagrams $D_1$ and $D_2$ (or, in other
words, only one special diagram $D=D_1 \sqcup D_2$, which is
disconnected sum of diagrams $D_1$ and $D_2$) is constructed  by
joining of following parts of arcs which belongs to $D$. The
boundary of the planar graph of $D$ consists of some edges of the
graph. Each edge is a part of an arc of $D$ (or the whole arc). We
choose a pair of edges in the boundary of $D$ which corresponds to a
pair of pre-critical arcs (or post-critical arcs) of $D=D_1 \sqcup
D_2$ (in first case we talk about {\em pre-critical sum} of special
diagrams, in the second case we talk about {\em post-critical sum})
and connect this pair, by using pair of lines which oriented in the
opposite manner. We call this  connected sum of special virtual
diagrams  {\em special}.

It is clear that a special connected sum of special diagrams is
special, too. Besides, each critical arc of $D$ remains critical in
the special connected sum for $D$.

One of  important characteristics of a special diagram is its
$T$-matrix. So, it is important to compare the $T$-matrix of a
special diagram $D$ and the $T$-matrix of its special sum
$\widetilde{D}$. Let classical crossings and long arcs of diagrams
$D$ and $\widetilde{D}$ be enumerated correspondingly, i.e. the
$i$-th long arc is emanating from the $i$-th classical crossing,
$i=1,...,n$, and identical classical crossings of $D$ and
$\widetilde{D}$ have the same numbers. Let $b_i$ and $b_j$ be a pair
of pre-critical (or post-critical) arcs of $D$ which take part in
the special sum. $\gamma^i\supset b_i$ and $\gamma^j\supset b_j$
denotes long arcs of $D$. Let $a_i\subset\gamma^i$ and
$a_j\subset\gamma^j$ be critical arcs of $D$.
$\widetilde{a_i}\subset\widetilde{\gamma^i}$, $i=1,...,n$, denote
critical and long arcs of $\widetilde{D}$.

So, if the special sum, which is constructed by  joining the arcs
$b_i$ and $b_j$, is pre-critical, then $\widetilde{a}_j=a_i$,
$\widetilde{a}_i=a_j$, consequently, $T(D)$ and $T(\widetilde{D})$
differ only by the transposition of the $i$-th and the $j$-th
columns. But if the special sum is post-critical then $a_i \subset
\widetilde{\gamma}^i$ and $a_j \subset \widetilde{\gamma}^j$, hence
$T(D)=T(\widetilde{D})$.

Assume the diagram $D$ is composed by two split pieces
$D_1$ and $D_2$. Then  $T(D)=\begin{pmatrix} T(D_1) & 0\\
0 & T(D_2)
 \end{pmatrix}$ for a convenient enumeration of classical crossings of $D$. Thus,
 the above argument yields the following equalities for the
 connected sum $\widetilde{D}=D_1+D_2$:
 $$det\, T(D_1+D_2)=\pm \,det\, T(D)=\pm\, det\,T(D_1)\, det\,
T(D_2)$$ $$per M(D_1+D_2)= per\,M(D)= per M(D_1)\,per M(D_2)$$

Consequently, the  $T$- and $M$-properties are invariant under the
connected sum operation (Note that the connected sum of $M$-diagrams
does not create cyclic crossings, since a set of critical arcs of a
connected sum are obtained by the union of the sets of critical arcs
 of the summands, hence,
the critical arc which is incident to a classical crossing three
times, can not appear unless it exists in one of the summands).

Thus, the connected sum operation described above allows to
construct new $M$-diagrams from the old ones. Which of them form new
$M$-links? To answer this question, it is useful to understand the
behaviour of $\beta$ and $\varepsilon$ under the connected sum
operation for $M$-diagram.

Note that for  $M$-diagrams, the connected sum preserves not only
critical arcs, but also the pairings: classical crossing
--- critical arc. Thus, the invariant
$\beta$ which is equal to the number of classical crossings paired
with overpassing critical arcs, satisfies the following equality:
$$\beta(D_1+D_2)=\beta(D_1)+\beta(D_2)$$

The same argument yields $x(D_1+D_2)=x(D_1)+x(D_2)$,
$y(D_1+D_2)=y(D_1)+y(D_2)$, where $x$ is the number of negative
classical crossings paired with critical arcs going over, and $y$ is
the number of classical crossings paired with incoming critical
arcs. Besides, $det\,M(D_1+D_2)=-det\,M(D_1)\,det\,M(D_2)$, if the
connected sum is taken with respect to pre-critical short arcs, and
$det\,M(D_1+D_2)=det\,M(D_1)\,det\,M(D_2)$, if it is taken with
respect to post-critical ones. This yields
$\varepsilon(D_1+D_2)=det\,M(D_1+D_2)\,(-1)^{x(D_1+D_2)+y(D_1+D_2)}=$
$\pm \,\varepsilon(D_1)\varepsilon(D_2)$, where the sign $-$ or $+$
is chosen according to the connected sum type: we take $-$ for a
pre-critical one and $+$ for the post-critical one.

The invariants $\varepsilon$ and $\beta$ and the connected sum
construction of $M$-diagram allow to construct the following four
series of $M$-knots.

\begin{figure}
\centering\includegraphics[scale = 0.8]{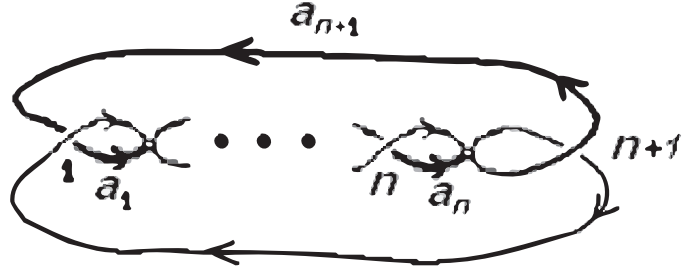}
  \caption{The $M$-knot $\Omega_n$ }
\label{fig5}
\end{figure}

\noindent 1) The $\Omega$-series are shown in Fig. \ref{fig5}.

We have
$M(\Omega_n)=\begin{pmatrix} 1 & 0& \ldots& 0& 1\\
                0  & 1 & 0 &\ldots& 0\\
                \hdotsfor{5}\\
                0&\ldots&0&1&0\\
                0&\hdotsfor{2}&0&1\\
                \end{pmatrix}$,
                $det\,M=1$, $x=1$, $y=0$ $\Rightarrow$
                $\varepsilon(\Omega_n)=det\,M\,(-1)^x\,(-1)^y=-1$.
                $\beta(\Omega_n)=1$.

\begin{figure}
\centering
  \includegraphics[scale = 0.7]{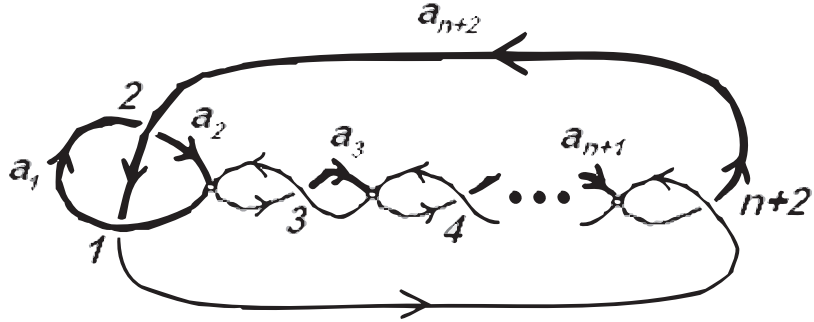}
  \caption{The $M$-knot $W_n$ }
\label{fig6}
\end{figure}

\noindent 2) The $W$-series are shown in Fig. \ref{fig6}.

We have
$M(W_n)=\begin{pmatrix} 1 & 0& \ldots& \ldots & 0& 1\\
                1  & 1 & 0 &\ldots& 0 & 1\\
                0 & 0 & 1 & 0 & \ldots & 0\\
                \vdots & \vdots &  \vdots& \ddots &\vdots&\vdots\\
                0&0&0&\ldots&1&0\\
                0& 0 & 0&\ldots&0&1\\
                \end{pmatrix}$, $det\,M=1$, the classical crossing
                $v_i$ is paired with the critical arc $a_i$,
                $i=1,...,n+2$, thus $\beta(W_n)=1$; $x=0$, $y=0$,
                consequently
                $\varepsilon(W_n)=det\,M\,(-1)^x\,(-1)^y=1$.

\begin{figure}
\centering
  \includegraphics[scale = 0.53]{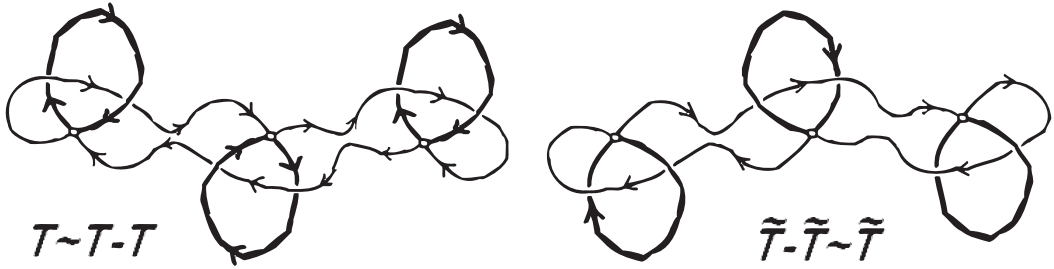}
  \caption{Examples of knots: $T$-series and their schemes}
\label{fig7}
\end{figure}

\noindent 3) The  $T$-series are obtained from $M$-knots $T$ and
$\widetilde{T}$ by using connected sums. This connected sum is
represented schematically by a graph. The vertices of the graph
correspond to virtual diagrams (summands, in our case $T$ and
$\widetilde{T}$), and the edges correspond to the ``tubes''
$S^{0}\times I$ for the connected sums; these edges are of two
types: pre-critical and post-critical, which is graphically
represented by a wave line or by a straight line, respectively, see
Fig. \ref{fig7}.

Clearly, there are virtual diagrams with schemes $\underbrace{T\sim
T-T\sim T-...-(\sim)T}_n$ and
$\underbrace{\widetilde{T}-\widetilde{T}\sim\widetilde{T}-\widetilde{T}\sim...-(\sim)\widetilde{T}}_n$.
Denote these diagrams by $T_n$ and ${\widetilde{T}}_n$,
respectively. According to the above property of connected sums of
$M$-diagrams, we have $\varepsilon(A\sim
B)=-\varepsilon(A)\,\varepsilon(B)$,
 $\varepsilon(A-B)=\varepsilon(A)\,\varepsilon(B)$. By
$\varepsilon(T)=1$, $\varepsilon(\widetilde{T})=-1$, this yields
$\varepsilon(T_n)=1$ $\Leftrightarrow$ $n{\equiv}_4\, 0$ or $1$, and
$\varepsilon(\widetilde{T}_n)=1$ $\Leftrightarrow$ $n{\equiv}_4\, 2$
or $3$, thus $T_n\nsim \widetilde{T}_n $ for all $n$.

\begin{figure}
\centering
  \includegraphics[scale = 0.53]{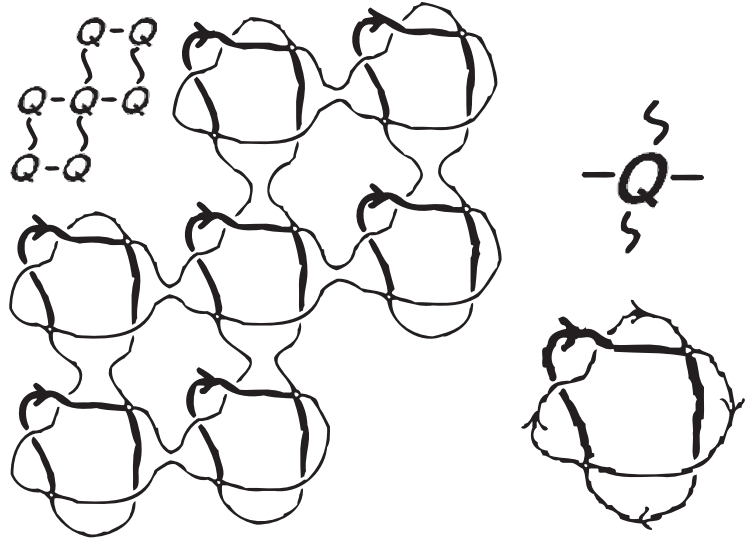}
  \caption{Constructing the  $Q$-series of knots}
\label{fig8}
\end{figure}

\noindent 4) The $Q$-series are constructed by using connected sums
of the diagram $Q$ shown in Fig. \ref{fig4}. By applying the
connected sum operation to $Q$, we can attach its parallel copies in
four directions, see Fig. \ref{fig8}.

As we see, in order to get a one-component link, it is not necessary
to take a tree for a scheme of a $Q$-knot.

To distinguish the knots from $Q$-series, let us use the
$\varepsilon$-invariant. Since $\varepsilon(Q)=1$, the value of
$\varepsilon$ on a connected sum of $r$ copies of $Q$ is equal to
$+1$ if and only if the number of wave lines in the scheme is even.
Thus, for $r\geqslant 2$ there are at least two knots in the
$Q$-series, for which the minimal virtual crossing number is equal
to $2r$, since in the scheme of connected sum we may change the
parity of the number of wave lines.

The information about these series is collected in the following
table.

{\footnotesize

\begin{tabular}{|p{3cm}|p{1cm}|p{1cm}|p{1cm}|p{2cm}|}
\hline
   & \quad $\Omega_n$ &  \ \ $W_n$ & $T_n,\widetilde{T}_n$ & $Q_n,n=2r$\\

\hline \begin{center}{\LARGE$\varepsilon$} \end{center}&
\begin{center}$-1$ \end{center}&
\begin{center}$1$ \end{center}& \begin{center}$\pm 1$ \end{center}&  \begin{center}$\pm 1$\end{center}\\

\hline \begin{center}{\LARGE $\beta$ }\end{center}&
\begin{center}$1$ \end{center}&
\begin{center}$1$ \end{center}& \begin{center}{\large $n$} \end{center}& \begin{center}{\Large $\frac{n}{2}$}\end{center}\\

\hline The number of classical crossings  &
\begin{center}{\large $n+1$} \end{center} & \begin{center}{\large  $n+2$}\end{center} &
\begin{center}{\large$2n$} \end{center}& \begin{center}{\large $\frac{3}{2}n$} \end{center} \\
\hline
\end{tabular}

}

Hence, for $n\geqslant 2$ the knots $\Omega_n, W_n, T_n,
\widetilde{T}_n$ are pairwise distinct, and for even $n\geqslant 4$
the knots $\Omega_n, W_n, T_n, \widetilde{T}_n, Q^{\prime}_n,
Q^{\prime\prime}_n$ are pairwise distinct, where $Q^{\prime}_n,
Q^{\prime\prime}_n$ are any two $Q$-knots with minimal virtual
crossing  number  $n$. It is easy to check that for $n=1$ the series
listed above consist of the following three knots: $T_1,
\widetilde{T}_1\sim\Omega_1, W_1$.

\section{Discussion and open questions.}

It is known that classical links embed into virtual links. More
precisely, if two classical diagrams are connected by a chain of
virtual Reidemeister moves, then they can be connected by a chain of
classical Reidemeister moves (see \cite{GPV}).

There is an analogous question between the relation of virtual links
generated by $M$-diagrams and all virtual links.

If one wants to classify some objects, say, virtual links, it would
be natural to split the object of classification into some
``layers'': say, classical links (those having no virtual crossings
at all), links with one virtual crossing, links with two virtual
crossing etc. Possibly, one may tackle this problem by using a
subset of links represented by $M$-diagrams.

Therefore, we list some open questions concerning $M$-property.

\noindent {\scshape Conjecture 1.} {\em $M$-property of proper
virtual diagrams is stable with respect to the virtual Reidemeister
moves
${\Omega}_1,{\Omega}_2,{\Omega}_3,{\Omega}^{\prime},{\Omega}^{\prime}_3$.}

\noindent {\scshape Conjecture 2.} {\em If virtual diagrams are
equivalent and satisfy $M$-property, then they are connected by
Reidemeister moves
${\Omega}_1,{\Omega}_2,{\Omega}_3,{\Omega}^{\prime},{\Omega}^{\prime}_3$.}

Moreover, it can be formulated more strong

\noindent {\scshape Conjecture 3.} {\em If virtual diagrams $D$ and
$D'$ are equivalent, $D$ satisfies $M$-property and $D'$ is minimal,
then $D'$ satisfies $M$-property.}

Finally, we address the following two questions.

1) {\em Is the number of nonequivalent $M$-diagrams which have the
same fixed number of virtual crossings, finite? }

2){\em Is it true that the number of classical crossing for a
one-component $M$-diagram grater than the number of virtual
crossings?}


\begin{thebibliography}{100}

\bibitem[BF]{BF} A.\,~Bartholemew, R.\,~Fenn (2003), Quaternionic
invariants of virtual knots and links, {\em Journal of Knot Theory
and Its Ramifications}, {\bf  17} (2),2008 pp. 231-251


\bibitem[DK]{DK} H.\,~Dye, L.H.\,~Kauffman, Virtual Crossing Number
and the arrow Polynomial, Arxiv.Math/GT. 0810:3858.


 \bibitem[FKM]{FKM}
R.\,A.~Fenn, L.\,H.~Kauffman, V.\,O.~Manturov, Virtual knots -
unsolved problems (2006), Fundamenta Mathematicae, N.\, 188, pp.\,
293--323.

 \bibitem[GPV]{GPV}
M.~Goussarov, M.~Polyak, and O.~Viro, Finite type invariants of
classical and virtual knots, Topology. 2000. V.\, 39, pp.\,
1045--1068.

 \bibitem[Jon]{Jones}
V.\,F.\,R.~Jones, A polynomial invariant for links via Neumann
algebras (1985), {\em Bull.\, Amer.\,  Math.\,  Soc.}, {\bf 129},
pp.\, 103--112.

 \bibitem[Ka1]{KaV}
L.\,H.~Kauffman, Virtual knot theory, Eur. J. Combinatorics. 1999.
V.\, 20, N.\, 7, pp.\, 662--690.

 \bibitem[Ka2]{KauffmanBracket}
L.\,H.~Kauffman, State Models and the Jones Polynomial (1987), {\em
Topology}, {\bf 26}, pp.\, 395--407.



\bibitem[KM]{KM} L.H.Kauffman, V.O.Manturov,
Virtual Knots and Links, Proceedings of the Steklov Inst, 2006, Vol.
252, N.1, pp. 104-121.

\bibitem[KR]{KR} L.H.Kauffman, D.Radford (2002), Bi-oriented quantum
algebras and a generalized Alexander polynomial for virtual links,
{\em AMS Contemp. Math.}, {\bf 318}, pp. 113-140.

 \bibitem[Kup]{Kup}
G.~Kuperberg, What is a Virtual Link? (2002), www.arXiv.org,
math-GT$\slash$ 0208039, {\em Algebraic and Geometric Topology},
2003, {\bf 3}, 5 87-591.

\bibitem[Ma1]{MaXi}
 V.O.Manturov, An Invariant 2-variable polynomial for virtual links (2002),
(Russian Math. Surveys), {\bf 57}, No.5, P.141-142.

\bibitem[Ma2]{Minim} V.O.Manturov, The Khovanov complex and minimal knot diagrams
(2006), Doklady Mathematics, {\bf 73}, N.1, P. 46-48.

 \bibitem[Ma3]{MyBook}
V.\,O.~Manturov, {\em Teoriya Uzlov} (Knot Theory), (Moscow-Izhevsk,
RCD), 2005 (512 pp).

\bibitem[Miy]{Miyazawa} Miyazawa, Y (2006), Magnetic Graphs and an
Invariant for virtual links, {\em J. of Knot Theory and Its
Ramifications}, {\bf 15} (10), pp.1319-1334.

\bibitem[Mur]{Mur}
K.\,~Murasugi, The Jones Polynomial and Classical Conjectures in
Knot Theory, {\em Topology}, {\bf 26}, 1987, pp. 187-194.

\bibitem[Saw]{Saw} J. Sawollek (2002), On Alexander-Conway Polynomials for
Virtual Knots and Links, {\em J. Knot Theory and Its Ramifications},
{\bf 12} (6), pp.767-779.

\bibitem[SW]{SW} D.Silver and S.Williams (2001), Alexander Groups
and Virtual Links, {\em J. of Knot Theory and Its Ramifications},
{\bf 10} (1), pp. 151-160.


\bibitem[Tep]{Tep} E.V.Teplyakov, On Roots of the $\xi$-polynomial (2005),
{\em Journal of Mathematical Sciences}, {\bf 128} (4), P. 3107-3113.

\bibitem[Thi]{Thi}
M.\,~Thistlethwaite, A spanning tree expansion for the Jones
polynomial, {\em Topology}, {\bf 26}, 1987, pp. 297-309.




\end{thebibliography}
 \end{document}